\UseRawInputEncoding
\documentclass[12pt]{article}

\usepackage{mathrsfs}
\usepackage{amssymb}
\usepackage{amsmath}
\usepackage{amsfonts}
\usepackage{amstext}
\usepackage{amsthm}
\usepackage{graphicx}
\usepackage{subfig}
\usepackage{color}
\usepackage{indentfirst,latexsym,bm}
\usepackage{mathrsfs}
\usepackage{cite}
\usepackage[all]{xy}
\usepackage{newlfont}

\usepackage{enumerate}
\numberwithin{equation}{section}

\numberwithin{equation}{section}

\begin{document}

\title{Interval bifurcation theorems for Fredholm operator and its application to an elliptic overdetermined problem
\thanks{Research supported by National Natural Science Foundation of China (No. 12371110, No. 12301133), the Postdoctoral Science Foundation of China (No. 2023M741441, No. 2024T170353) and Jiangsu Education Department (No. 23KJB110007).}}

\author{
Guowei Dai
\thanks{ School of Mathematical Sciences, Dalian University of Technology, Dalian, 116024, P.R. China. E-mail: daiguowei@dlut.edu.cn.} \,,
Yong Zhang
\thanks{Corresponding author. School of Mathematical Sciences, Jiangsu University, Zhenjiang, 212013, P.R. China. E-mail: zhangyong@ujs.edu.cn}
}
\date{}
\maketitle

\renewcommand{\abstractname}{Abstract}

\renewcommand{\abstractname}{Abstract}

\begin{abstract}
We establish local interval bifurcation theorem and global interval bifurcation theorem for Fredholm operator with index $0$ via $0$-group.
As one of applications, we investigate the existence of a family of nontrivial domains $\Omega_{\rho}\subset \mathbb{R}^N$ ($N=2,3$ or $4$), bifurcating from a small ball, such that the problem
\begin{equation}
-\Delta u=u-\left(u^+\right)^3\,\, \text{in}\,\,\Omega_{\rho}, \,\, u=0,\,\,\partial_\nu u=\text{const}\,\,\text{on}\,\,\partial\Omega_{\rho} \nonumber
\end{equation}
has a sign-changing bounded solution. Compared with the recent result \cite[Theorem 2.1]{Ruiz}, here we obtain a family of domains $\Omega_{\rho}$ instead of a sequence of domains.
\end{abstract}


\emph{AMS Subjection Classification(2020):} 35N05; 47A53; 47A55; 47J15


\section{Introduction}
\quad\,
Let $X$ and $Y$ be real Banach spaces with $X\subseteq Y$. We investigate the structure of the set of nontrivial solutions for the following nonlinear operator equation
\begin{equation}\label{operater}
F(\lambda,u)=0,\,\, (\lambda,u)\in \mathbb{R}\times X,
\end{equation}
where $F:\mathbb{R}\times X\rightarrow Y$
is continuous with $F(\lambda,0)=0$ for $\lambda\in \mathbb{R}$ and $C^1$ with respect to $u$ at $u=0$.
Kielh\"{o}fer \cite{Kielhofer} established the bifurcation theory to (\ref{operater}) via odd crossing number.
We next recall the conception of odd crossing number.

Let $0$ be an isolated eigenvalue of algebraic multiplicity $m$ of $D_uF\left(\mu,0\right)$ for some $\mu\in \mathbb{R}$.
It is well known (see \cite{Kato} or \cite{Kielhofer}) that
the number $m$ is an invariant, i.e., the dimension of eigenspace is invariant under perturbation near $\mu$.
The set of all perturbed eigenvalues near $0$ is called \emph{$0$-group}.
Further, define $\Sigma(\lambda)=1$ if there are no negative real eigenvalues in the $0$-group of $D_uF(\lambda,0)$ and
\begin{equation}
\Sigma(\lambda)=(-1)^{m_1+m_2+\dots+m_k}\nonumber
\end{equation}
if $\mu_1,\mu_2,\dots,\mu_k$ are all negative real eigenvalues in the $0$-group having algebraic multiplicities $m_1,m_2,\dots,m_k$, respectively.
From now on, for simplicity, $\Sigma(\lambda)$ is called \emph{$0$-group index} of $\lambda$.
If $D_uF(\lambda,0)$ is regular for $\lambda\in\left(\mu-\delta, \mu\right) \cup\left(\mu, \mu+\delta\right)$ and if $\Sigma(\lambda)$ changes at $\lambda=\mu$,
then $D_uF(\lambda,0)$ has an {odd crossing number} at $\lambda=\mu$.

If $D_u F(\lambda,0):X\longrightarrow Y$ is a Fredholm map of index $0$ for $\lambda$ near $\mu$ and $D_u F\left(\lambda,0\right)$ has an odd crossing number at $\lambda=\mu$,
Kielh\"{o}fer \cite[Theorem II.4.4]{Kielhofer} proved that $\left(\mu, 0\right)$ is a bifurcation point for $F(\lambda,u)=0$.
Furthermore, if $F$ is $C^2$, $D_u F(\lambda,u):X\longrightarrow Y$ satisfies the admissible condition of \cite[Definition II.5.1]{Kielhofer} for all $\lambda\in \mathbb{R}$ and $F$ is proper on every closed bounded subset of $\mathbb{R}\times X$, he \cite[Theorem II.5.8]{Kielhofer} also
proved that there exists a global continua $\mathscr{C}_\mu$ bifurcating from $\left(\mu, 0\right)$, such that $\mathscr{C}_{\mu}$ satisfies Rabinowitz-type
global alternative.

The \emph{first aim} of this work is to extend Kielh\"{o}fer's local point bifurcation result into interval bifurcation result as follows.
\\ \\
\textbf{Theorem 1.1.} \emph{Let $X$ and $Y$ be real Banach spaces with $X\subseteq Y$ and $F:\mathbb{R}\times X\longrightarrow Y$ be a $C^1$ map with $F(\cdot,0)=0$. Suppose that}

1. \emph{For some $a,b\in \mathbb{R}$ with $a\leq b$, $D_{u}F(\lambda, 0)$ is a Fredholm operator with index $0$ near $[a,b]$},

2. \emph{$D_uF(\lambda,0)$ is regular for $\lambda\in\left(a-\varepsilon, a\right) \cup\left(b, b+\varepsilon\right)$ for any $\varepsilon>0$ small enough},

3. \emph{$\Sigma (a-\varepsilon)\neq \Sigma (b+\varepsilon)$ for any $\varepsilon>0$ small enough}.

\noindent \emph{Then $[a,b]\times \{0\}$ is a bifurcation interval of $F(\lambda,u)=0$ in the sense: every neighborhood of
$[a,b]\times \{0\}$ contains solutions of $F(\lambda,u)=0$ with $(\lambda,u)\in \mathbb{R}\times X\setminus\{0\}$}.
\\

When $a=b$, Theorem 1.1 degenerates to Kielh\"{o}fer's local point bifurcation result \cite[Theorem II.4.4]{Kielhofer}. The \emph{second aim} is to extend Kielh\"{o}fer's global point bifurcation result into interval bifurcation by degree theory for $C^1$ Fredholm mappings of index $0$.
Let $\mathcal{S}$ be the closure of the set of nontrivial solutions of equation (\ref{operater}).
\\ \\
\textbf{Theorem 1.2.} \emph{Let $X$ and $Y$ be real Banach spaces with $X\subseteq Y$ and $F:\mathbb{R}\times X\longrightarrow Y$ be a $C^1$ map with $F(\cdot,0)=0$. Suppose that}

1. \emph{$D_{u}F(\lambda, u)$ is a Fredholm operator with index $0$ for all $(\lambda,u)\in \mathbb{R}\times X$},

2. \emph{$D_uF(\lambda,0)$ is regular for $\lambda\in\left(a-\varepsilon, a\right) \cup\left(b, b+\varepsilon\right)$ for any $\varepsilon>0$ small enough},

3. \emph{$\Sigma (a-\varepsilon)\neq \Sigma (b+\varepsilon)$ for any $\varepsilon>0$ small enough},

4. \emph{the restriction of $F$ to closed bounded subsets of $\mathbb{R}\times X$ is compact}.

\noindent \emph{Then $\mathcal{S}$
possesses a maximal component $\mathscr{C}_\mu$ emanating from $[a,b]\times \{0\}$, such that either $\mathscr{C}_{\mu}$ is unbounded or contains some $\left(\lambda_*,0\right)$ with $\lambda_*\in \mathbb{R}\setminus[a,b]$}.
\\

When $a=b$, the conclusion of Theorem 1.2 degenerates to Kielh\"{o}fer's global point bifurcation result \cite[Theorem II.5.8]{Kielhofer}.
In fact, even in the case of $a=b$, our conditions are weaker than those in \cite[Theorem II.5.8]{Kielhofer}.
Because we do not require $F$ being $C^2$ and the spectrum condition \cite[Definition II.5.1]{Kielhofer}.

As one of applications of Theorem 1.2, we will discuss Serrin's overdetermined problem on nontrivial bounded domain in $\mathbb{R}^N$.
It is known that the following overdetermined elliptic problem
\begin{equation}\label{ffunction}
\left\{
\begin{array}{ll}
\Delta u+f(u)=0\,\, &\text{in}\,\, \Omega,\\
u=0 &\text{on}\,\, \partial \Omega,\\
\frac{\partial u}{\partial \nu}=\text{const} &\text{on}\,\, \partial \Omega
\end{array}
\right.
\end{equation}
has attracted a lot of attention, where $\Omega$ is a given regular domain, $f$ is a Lipschitz function and $\nu$ is the unit outward normal about $\partial \Omega$. In fact, it is natural to consider this class of problem as they would appear in many different physical phenomena \cite{Sirakov}.
In Euclidean space $\mathbb{R}^N$, J. Serrin proved that, if $\Omega$ is a bounded domain of class $C^2$ such that (\ref{ffunction}) admits a solution, then
 $\Omega$ is a ball.
In 1997, under some hypothesis on the nonlinearity $f(u)$ and on the behaviour of the epigraph at infinity, Berestycki, Caffarelli and Nirenberg \cite{BCN} proved that if problem (\ref{ffunction}) admits a smooth and bounded solution, then the epigraph $\Omega\subseteq \mathbb{R}^N$ is a half-space. In addition, they also proposed the following famous conjecture in \cite{BCN}.
\\ \\
\textbf{BCN Conjecture:} If $\Omega$ is a smooth domain and $\mathbb{R}^N\setminus \overline{\Omega}$ is connected such that problem (\ref{ffunction}) exists a bounded solution, then $\Omega$ is either a ball, a half-space, a generalized cylinder $B^k\times \mathbb{R}^{N-k}$ where $B^k$ is a ball in $\mathbb{R}^k$, or the complement of one of them.\\

\noindent These results were later extended by A. Farina and E. Valdinoci \cite{Farina}, where they proved that a globally Lipschitz smooth epigraph $\Omega$ of $\mathbb{R}^N$ ($N=2$, $3$) admitting a positive solution $u$ to problem (\ref{ffunction}) with $C^1$ bistable nonlinearity, must be a half-space and $u$ must be a function of only one variable.
With some more weaker conditions on $f$ and $u$ for $N=2$, A. Ros and P. Sicbaldi \cite{Ros} proved that the BCN Conjecture is valid on $C^{1,\alpha}$ domain whose boundary is unbounded and connected (also see \cite{RRS1}). We also like to point that Serrin's result has been generalized to exterior domains (see \cite{Aftalion, Reichel, Sirakov, Sirakov1}) and ring-shaped domain (see \cite{Reichel0, Willms}).

On the other hand, we find that constructing the counterexamples to BCN Conjecture has attracted much attention recently. The breakthrough in this direction in the Euclid space $\mathbb{R}^N$ is due to P. Sicbaldi \cite{Sicbaldi}, where he constructed the first counterexample to BCN Conjecture and proved that the cylinder $B_1\times \mathbb{R}$ can be perturbed to an unbounded domain, so that problem (\ref{ffunction}) with $f(s)=\lambda s$ has a bounded positive
solution, where $B_1$ is the unit ball of $\mathbb{R}^N$ with $N\geq2$ centered at the origin.
F. Schlenk and P. Sicbaldi \cite{Schlenk} further proved that these new extremal domains belong to a smooth bifurcation family of domains for $N\geq1$. When $N\geq 9$, M. Del Pino, F. Pacard and J. Wei \cite{Del} established an epigraph domain which is not a half-space, such that the problem (\ref{ffunction}) with Allen-Cahn type nonlinearities is solvable. In \cite{Fall}, the authors also gave a counter-example to the BCN's
conjecture with $f(s)\equiv1$ and A. Ros, D. Ruiz and P. Sicbaldi \cite{RRS} obtained an important counter-example to this
conjecture on exterior domain which is not the complement of a ball.

All the previous discussion is on the positive solution of (\ref{ffunction}). In fact, the sign-changing solutions are also significant in the applications.
For instance, Ruiz \cite{Ruiz} recently found a sign-changing solution to the problem for a special function $f(u)$ in a bounded domain $\Omega$ different from a ball. Considering the following problem
\begin{equation}\label{canshu}
\left\{
\begin{array}{ll}
-\rho\Delta u= u-\left(u^{+}\right)^3 \,\, &\text{in}\,\, \Omega,\\
u=0 &\text{on}\,\, \partial \Omega,\\
\frac{\partial u}{\partial \nu}=\text{const} &\text{on}\,\, \partial \Omega,
\end{array}
\right.
\end{equation}
he used the Krasnoselskii bifurcation theorem to show the existence of a sequence of nontrivial sign-changing solutions to (\ref{canshu}). We will strengthen the result and obtain a global branch of nontrivial solutions to (\ref{canshu}) by using Theorem 1.2.

To state our main result, let us first introduce some notations. Denote by $B(R)\subset\mathbb{R}^N$ the ball centered at the origin of radius $R$ and $B=B(1)$. We write $G$ to denote a symmetry group $G\subset O(N)$ and a function $f$ is $G$-symmetric if $f\circ g=f$ for any $g\in G$.
Denote by $\lambda_i$ the eigenvalues of the Laplacian operator on the unit ball $B$ with Dirichlet boundary conditions for $G$-symmetric functions, counted with multiplicity. We also denote by $\bar{\lambda}_i$, the Dirichlet eigenvalues associated to radial eigenfunctions. Of course we have that
\begin{equation}
0<\lambda_1<\lambda_2\leq \lambda_3....,\quad 0<\bar{\lambda}_1<\bar{\lambda}_2<\bar{\lambda}_3....\nonumber
\end{equation}
and
\begin{equation}
\lambda_1=\bar{\lambda}_1, \lambda_i\leq \bar{\lambda}_i~~\text{if}~~i>1.
\nonumber
\end{equation}
In the similar way, we denote by $\sigma_i$ and $\bar{\sigma}_i$ the eigenvalues of the Laplacian operator on the unit ball with Neumann boundary conditions for $G$-symmetric and radial functions, respectively. There holds that
\begin{equation}
0=\sigma_0<\sigma_1\leq \sigma_2\leq \sigma_3....,\quad 0=\bar{\sigma}_0<\bar{\sigma}_1<\bar{\sigma}_2< \bar{\sigma}_3....
\nonumber
\end{equation}
and
\begin{equation}
0=\bar{\sigma}_0<\bar{\lambda}_1<\bar{\sigma}_1<\bar{\lambda}_2<\bar{\sigma}_2<\bar{\lambda}_3....
\nonumber
\end{equation}
From now on, we fix a symmetry group $G$ of $\mathbb{S}^{N-1}$ with the following property:\\

(G) Denote by $\sigma$ the first eigenvalue $\sigma_{k}$ with $\sigma\neq\bar{\sigma}_j$ for all $j\in \mathbb{N}$ and assume that $\sigma>\bar{\lambda}_2$ and has one multiplicity.\\

In fact, there exist the symmetry group $G$ satisfying (G) (see the Appendix in \cite{Ruiz}). Let us recall that the eigenvalues of the Laplace-Beltrami operator on $\mathbb{S}^{N-1}$ have the expression $i(i+N-1)$ for $i\in \mathbb{N}$. We denote by $\mu_{k}$ the eigenvalues of the Laplace-Beltrami operator on $\mathbb{S}^{N-1}$ for $G$-symmetric eigenfunctions. Obviously, $\mu_{i_k}=i(i+N-2)$ for some $i_k\in \mathbb{N}$. The associated $G$-symmetric eigenfunctions are denoted by $\zeta_{k}$ and are normalized such that
$\int_{\mathbb{S}^{N-1}}\zeta_{k}^2(\theta)\,\text{d}\theta=1$.
Define the spaces of H\"{o}lder continuous functions
\begin{equation}
C^{k,\alpha}_{G}(B(R))=\left\{ u\in C^{k,\alpha}(B(R)): u~\text{is}~G-\text{symmetric} \right\},
\nonumber
\end{equation}
\begin{equation}
C^{k,\alpha}_{G,0}(B(R))=\left\{ u\in C^{k,\alpha}_{G}(B(R)): u=0~\text{on}~\partial B \right\},
\nonumber
\end{equation}
\begin{equation}
C^{k,\alpha}_{G}(\mathbb{S}^{N-1})=\left\{ u\in C^{k,\alpha}(\mathbb{S}^{N-1}): u~\text{is}~G-\text{symmetric} \right\}
\nonumber
\end{equation}
and
\begin{equation}
C^{k,\alpha}_{G,m}(\mathbb{S}^{N-1})=\left\{ u\in C^{k,\alpha}_{G}(\mathbb{S}^{N-1}): \int_{\mathbb{S}^{N-1}}u\,\text{d}\theta=0 \right\}.
\nonumber
\end{equation}
At last, let us introduce two useful spaces of functions that are $L^2$-orthogonal to all radial functions by
\begin{equation}
E=\left\{ \psi\in H^{1}_{G}(B): \int_{B}\psi(x)g(x)\,\text{d}x=0~~\forall g\in L^2_{r}(B) \right\}
\nonumber
\end{equation}
and
\begin{equation}
E_0=\left\{ \psi\in H^{1}_{0,G}(B): \int_{B}\psi(x)g(x)\,\text{d}x=0~~\forall g\in L^2_{r}(B) \right\}.
\nonumber
\end{equation}
Then, by Theorem 1.2 we can obtain a family of domains $\Omega_{\rho}\subset \mathbb{R}^N$ bifurcating from $B$ for which the problem (\ref{canshu})
has a sign-changing bounded solution.
\\ \\
\textbf{Theorem 1.3.} \emph{Let $G\subset O(N)$ be a symmetric group satisfying (G)}.
\emph{Then there exists a nontrivial branch $\mathcal{C}$ which emanates from some interval $\left[\rho_1,\rho_2\right]\times \{0\}$ in $\mathbb{R}\times C_{G,m}^{2,\alpha}\left(\mathbb{S}^{N-1}\right)$ with $0<\rho_1< \rho_2$ such that for each $(\rho,v)\in  \mathcal{C}\setminus\left(\left[\rho_1,\rho_2\right]\times \{0\}\right)$, the problem}
\begin{equation}
\left\{
\begin{array}{lll}
-\rho\Delta u=u-\left(u^+\right)^3 \,\, &\text{in}\,\, \Omega_{\rho},\\
u=0 &\text{on}\,\, \partial\Omega_{\rho},\\
\frac{\partial u}{\partial \nu}=\text{const} &\text{on}\,\, \partial\Omega_{\rho}
\end{array}
\right.\nonumber
\end{equation}
\emph{admits a sign-changing solution $u\in C_{G}^{2,\alpha}\left(\Omega_{\rho}\right)\cap H_{0,G}^{1}\left(\Omega_{\rho}\right)$, where $\Omega_{\rho}=\left\{ x\in \mathbb{R}^{N}: |x|<1+v\right\}$.}
\\

The outline of the rest of this article is as follows. In Section 2, we give the proofs of Theorems 1.1--1.2. In Section 3, we mainly recall some known results on Dirichlet problem obtained in \cite{Ruiz}. In Section 4, we investigate the properties of eigenvalues $\sigma_{i_1}\left(H_{\rho}\right)$ for the linearized operator $H_{\rho}$. Finally, in Section 5 we provide the proof of Theorem 1.3.

\section{Proofs of Theorems 1.1--1.2}

\quad\,
Let $\Phi_0(X,Y)$, $GL(X,Y)$ and $K(X,Y)$ be the set of linear Fredholm operators of index $0$, linear invertible operators and linear compact operators, respectively.
Given a continuous path $\alpha:[a,b]\longrightarrow \Phi_0(X,Y)$, it has been proved in \cite{Fitzpatrick0} that there exists a continuous path $\beta:[a,b]\longrightarrow GL(X,Y)$ such that
\begin{equation}
\beta(\lambda)\alpha(\lambda)-I\in K(X,X).\nonumber
\end{equation}
The path $\beta$ is called a \emph{parametrix} of $\alpha$.
Set $I-\beta(\lambda)\alpha(\lambda):=P(\lambda)$.
If $\alpha(a)$ and $\alpha(b)$ are isomorphisms, the \emph{parity} (see \cite{FitzpatrickP, FitzpatrickPejsachowicz, FitzpatrickRabier}) of $\alpha$ on $[a,b]$ is defined by
\begin{equation}
\sigma(\alpha,[a,b])=\deg(\beta(a)\alpha(a))\deg(\beta(b)\alpha(b)),\nonumber
\end{equation}
where $\deg(\beta(a)\alpha(a))$ is the Leray-Schauder degree of $(I-P(a))(x)=0$ with respect to any open set containing $0$.
\\ \\
\textbf{Proof of Theorem 1.1.} For any $\varepsilon>0$ small enough, we have that
\begin{equation}
\Sigma(a-\varepsilon)=(-1)^{m_1+m_2+\dots+m_k}\nonumber
\end{equation}
if $\mu_1,\mu_2,\dots,\mu_k$ are all negative real eigenvalues in the $0$-group having algebraic multiplicities $m_1,m_2,\dots,m_k$, respectively.
Similarly, we also have that
\begin{equation}
\Sigma(b+\varepsilon)=(-1)^{m_1'+m_2'+\dots+m_l'}\nonumber
\end{equation}
if $\mu_1',\mu_2',\dots,\mu_l'$ are all negative real eigenvalues in the $0$-group having algebraic multiplicities $m_1',m_2',\dots,m_l'$, respectively.
Set $m_1+m_2+\dots+m_k:=m_*$ and $m_1'+m_2'+\dots+m_l':=m^*$.
Since $\Sigma (a-\varepsilon)\neq \Sigma (b+\varepsilon)$, $m^*-m_*$ is an odd number.

We know that
\begin{equation}
\sigma\left(D_{u}F(\lambda, 0),[a-\varepsilon,b+\varepsilon]\right)=\deg\left(\beta(a-\varepsilon)D_{u}F(a-\varepsilon, 0)\right)\deg\left(\beta(b+\varepsilon)D_{u}F(b+\varepsilon, 0)\right),\nonumber
\end{equation}
where $\beta$ is a parametrix of $D_{u}F(\lambda, 0)$.
Since $\beta(\mu\pm\varepsilon)D_{u}F(\mu\pm\varepsilon, 0)$ is linear invertible, $I-P(\mu\pm\varepsilon)$ is linear invertible and completely continuous vector field, where $P(\lambda)=I-\beta(\lambda)D_{u}F(\lambda, 0)$. Using the index formula for isolated zeros \cite[Formula II.2.11]{Kielhofer}, we obtain that
\begin{equation}
\deg\left(\beta(a-\varepsilon)D_{u}F(a-\varepsilon, 0)\right)=(-1)^{m_*}\nonumber
\end{equation}
and
\begin{equation}
\deg\left(\beta(b+\varepsilon)D_{u}F(b+\varepsilon, 0)\right)=(-1)^{m^*}.\nonumber
\end{equation}
It follows that
\begin{equation}
\sigma\left(D_{u}F(\lambda, 0),[a-\varepsilon,b+\varepsilon]\right)=(-1)^{m_*+m^*}=(-1)^{m^*-m_*+2m_*}=(-1)^{m^*-m_*}=-1.\nonumber
\end{equation}
Applying \cite[Theorem 1]{Fitzpatrick}, we obtain that $[a,b]\times \{0\}$ is a bifurcation interval of $F(\lambda,u)=0$.
\qed\\

Although Theorem 1.1 was obtained through the application of \cite[Theorem 1]{Fitzpatrick}, odd crossing number may be easier to calculate than parity in some practical applications.
So, Theorem 1.1 may be easier to verify on some specific problems.
\\ \\
\textbf{Proof of Theorem 1.2.} From the argument of Theorem 1.1 we have known that
\begin{equation}
\sigma\left(D_{u}F(\lambda, 0),[a-\varepsilon,b+\varepsilon]\right)=-1.\nonumber
\end{equation}
The desired global interval conclusion can be obtained by applying \cite[Theorem 6.1]{Pejsachowicz}.\qed
\\

If $F$ is not globally defined, it is not difficult to get the following result.
\\ \\
\textbf{Corollary 2.1.} \emph{Assume that $\mathcal {O}$ is an open subset of $\mathbb{R}\times X$ and ${F}$ is defined on $\mathcal {O}$. Under the assumptions of Theorem 1.2, either $\mathscr{C}_{\mu}$ satisfies the alternatives of Theorem 1.2 or meets $\partial\mathcal {O}$.}\\

We end this section by providing a criterion for judging properness of $F$.
\\ \\
\textbf{Proposition 2.1.} \emph{Let $X$ and $Y$ be real Banach spaces and $F:\mathbb{R}\times X\longrightarrow Y$ be a $C^1$ map such that
$D_{u}F(\lambda, 0)$ is a Fredholm operator with index $0$ for all $\lambda\in \mathbb{R}$. If the embedding of $X\hookrightarrow Y$ is compact, the restriction of $F$ to closed bounded subsets of $\mathbb{R}\times X$ is proper}.
\\ \\
\textbf{Proof.} Assume that $u_n\rightharpoonup u_0$ and $\lambda_n\rightharpoonup \lambda_0$ with $F\left(\lambda_n,u_n\right)=y_n\rightarrow y_0$ in $Y$.
It is enough to show that $u_n\rightarrow u_0$ in $X$.
Since the embedding of $X\hookrightarrow Y$ is compact, we see that $u_n\rightarrow u_0$ in $Y$.
Thus, we have that $\lim_{n\rightarrow+\infty}F\left(\lambda_n,u_n\right)=F\left(\lambda_0,u_0\right)=y_0$.
Letting $H(\lambda,u)=F(\lambda,u)-D_u F(\lambda,0)u$, then we see that $H(\lambda,u)\in Y$.

Let $\beta:\mathbb{R}\rightarrow GL(Y,X)$ be a parametrix for $D_u F(\lambda,0)$.
Then, we find that $F(\lambda,u)=y$ if and only if $\beta(\lambda)F(\lambda,u)=\beta(\lambda)y$ for $(\lambda,u)\in \mathbb{R}\times X$ and $y\in Y$.
Let $L(\lambda)=\beta(\lambda)D_u F(\lambda,0)-I$.
It follows that
\begin{equation}
u=\beta(\lambda)F(\lambda,u)-L(\lambda)u-\beta(\lambda)H(\lambda,u)\nonumber
\end{equation}
for $(\lambda,u)\in \mathbb{R}\times X$.
In particular, there is
\begin{eqnarray}\label{proper1}
u_n&=&\beta\left(\lambda_n\right)F\left(\lambda_n,u_n\right)-L\left(\lambda_n\right)u_n-\beta\left(\lambda_n\right)H\left(\lambda_n,u_n\right)\nonumber\\
& =&\beta\left(\lambda_n\right)y_n-L\left(\lambda_n\right)u_n-\beta\left(\lambda_n\right)H\left(\lambda_n,u_n\right).
\end{eqnarray}
Note that
\begin{eqnarray}\label{proper2}
H\left(\lambda_n,u_n\right)&=&y_n-D_u F\left(\lambda_n,0\right)u_n\nonumber\\
&=&y_n+\left(D_u F\left(\lambda,0\right)-D_u F\left(\lambda_n,0\right)\right)u_n-D_u F\left(\lambda,0\right)u_n.
\end{eqnarray}
Since $F$ is $C^1$ and $u_n\rightharpoonup u$ in $X$, letting $n\rightarrow+\infty$ on the both sides of (\ref{proper2}), we have that
\begin{equation}
\lim_{n\rightarrow+\infty}H\left(\lambda_n,u_n\right)=y-D_u F\left(\lambda,0\right)u=H(\lambda,u).\nonumber
\end{equation}
Obviously, $L(\lambda)$ is continuous.
Combining the compactness of $L(\lambda)$ with the continuity of $\beta(\lambda)$ and $L(\lambda)$, we conclude from (\ref{proper1}) that $u_n$ converges to $u$ in $X$.
Therefore, the restriction of $F$ to closed bounded subsets of $\mathbb{R}\times X$ is proper.\qed

\section{Some known results on the Dirichlet problem}

\quad\, 
As shown in the introduction, we will prove the Theorem 1.3 by using a bifurcation argument. Thus it is necessary to consider the following Dirichlet problem
\begin{equation}\label{equationB}
\left\{
\begin{array}{ll}
-\rho\Delta u=u-\left(u^+\right)^3\,\, &\text{in}\,\,  B,\\
u=0 &\text{on}\,\, \partial B.
\end{array}
\right.
\end{equation}
It will be convenient to make a change of scale and consider the equivalent problem
\begin{equation}\label{equationBR}
\left\{
\begin{array}{ll}
-\Delta v=v-\left(v^+\right)^3\,\, &\text{in}\,\,  B(R),\\
v=0 &\text{on}\,\, \partial B(R)
\end{array}
\right.\nonumber
\end{equation}
for some suitable $R=\rho^{-\frac{1}{2}}$. If $v$ is the radial solution, we have that
\begin{equation}\label{eigenvalueonbal001}
\left\{
\begin{array}{lll}
-v''(r)-(N-1)\frac{v'(r)}{r}=v(r)-(v(r)^+)^3,\,\,\, r\in[0,R),\\
v(R)=0,\\
v'(R)=1.
\end{array}
\right.\nonumber
\end{equation}
Based on \cite[Proposition 3.1]{Ruiz}, the existence and asymptotic behavior of a radial sign-changing solution $u_{\rho}$ to (\ref{equationB}) are obtained as follows.
\\ \\
\textbf{Proposition 3.1.} \emph{For any $\rho\in\left(0, \bar{\lambda}^{-1}_{2}\right)$, there exists a radial sign-changing solution $u_{\rho}\in C^{4,\alpha}(\mathbb{R})$ to the problem $(\ref{equationB})$ satisfying that}

\emph{(1) The function $u_{\rho}(r)$ has a unique zero at a point $p_{\rho}\in (0,1)$ and}
\begin{equation}
\left\{
\begin{array}{ll}
u_{\rho}(r)>0\,\, &\text{for}\,\,  r\in\left[0,p_{\rho}\right),\\
u_{\rho}(r)<0\,\, &\text{for}\,\,  r\in\left(p_{\rho},1\right)
\end{array}
\right.\nonumber
\end{equation}
with $u'_{\rho}(1)>0$;

\emph{(2) The radial solution of $(\ref{equationB})$ satisfying (1) is unique;}

\emph{(3) If $\rho\rightarrow\bar{\lambda}^{-1}_{2}$, then $u_{\rho}\rightarrow 0$ in $C^4$ sense;}

\emph{(4) If $\rho\rightarrow0$, then $u_{\rho}\rightarrow 1$ in compact sets of $B$ in $C^4$ sense. Moreover, if we make the change of variables and define
\begin{equation}
\tilde{v}_{R}(r)=v_{R}(r-p_{R}), \quad p_{R}=\rho^{-\frac{1}{2}}p_{\rho},\nonumber
\end{equation}
we have that $\tilde{v}_{R}\rightarrow\tilde{v}_{0}$ in $C^4$ sense in compact sets of $(-\infty, \pi)$, where}
\begin{equation}
\tilde{v}_{0}=\left\{
\begin{array}{ll}
-\tanh\left(\frac{r}{\sqrt{2}}\right)\,\, &\text{for}\,\,  r\leq 0,\\
-\frac{1}{\sqrt{2}}\sin(r)\,\, &\text{for}\,\,  r\in(0,\pi].
\end{array}  \nonumber
\right.
\end{equation}
~\\
\indent If we use $\dot{u}_{\rho}$ to denote the derivative of $u_{\rho}$ with respect to $\rho$, it's easy to deduce that
\begin{equation}
\dot{u}_{\rho}(r)=-u'_{\rho}(r)\frac{r}{2R}\rho^{-\frac{3}{2}}=-u'_{\rho}(r)\frac{r}{2\rho},\nonumber
\end{equation}
which indicates that $u_{\rho}$ is differential with respect to $\rho$.

In order to solve (\ref{equationB}) on a perturbation domain, we consider its linearized eigenvalue problem near the radial solution $u_{\rho}$ under Dirichlet boundary conditions. We will consider both the radially symmetric and the $G$-symmetric case. Denote the linearized operator $L: H^{1}_{0,G}\rightarrow \left(H^{1}_{0,G}(B)\right)^{-1}$ by
\begin{equation}\label{e2.5}
L:=-\rho \Delta-1+3\left(u_{\rho}^{+}\right)^2.\nonumber
\end{equation}
Now consider the eigenvalue problems
\begin{equation}
L(\psi)=\mu_{k}\psi, \quad \psi\in H^{1}_{0,G}(B)\nonumber
\end{equation}
and
\begin{equation}
L(\psi)=\bar{\mu}_{k}\psi, \quad \psi\in H^{1}_{0,r}(B)
\nonumber
\end{equation}
with $\psi_{k}, \bar{\psi}_{k}$ being the associated eigenfunctions. It's known that $\mu_{1}=\bar{\mu}_{1}$ is simple and $\psi_{1}$ is positive. We also define the quadratic form $Q_{D}: H^{1}_{0,G}(B)\rightarrow \mathbb{R}$ by
\begin{equation}\label{e2.6}
Q_{D}(\phi)=\int_{B}\rho|\nabla\phi|^2-\phi^2+3\left(u_{\rho}^{+}\right)\phi^2\,\text{d}x\nonumber
\end{equation}
and let $\bar{Q}_{D}=Q_{D}|_{H^{1}_{0,r}(B)}$.
From Proposition 3.2 and Remark 3.3 in \cite{Ruiz}, we know that
\\ \\
\textbf{Proposition 3.2.} \emph{$\bar{\mu}_{1}<0<\bar{\mu}_{2}$ and the map $U: \left(0,\bar{\lambda}_{2}^{-1}\right)\rightarrow C^4(B)$, $U(\rho)=u_{\rho}$ is $C^1$. There exist $\lambda_{0}\in\left(0,\bar{\lambda}_{2}^{-1}\right)$ such that for any $\rho\in \left(\rho_{0},\bar{\lambda}_{2}^{-1}\right)$, then $\mu_{2}>0$, that is to say that the linearized operator $L$, acting on $H^{1}_{0,G}(B)$, is nondegenerate. Furthermore, the quadratic form $Q_{D}$ is positive in the space $E_{0}$.}
\\

We recall the perturbation result due to \cite[Proposition 4.1]{Ruiz}.
\\ \\
\textbf{Proposition 3.3.} \emph{Assume that $\rho\in \left(\rho_{0},\bar{\lambda}_{2}^{-1}\right)$, there exists a neighborhood $\mathcal{U}$ of $0$ in $C_{G,m}^{2,\alpha}\left(\mathbb{S}^{N-1}\right)$ such that for any $v\in \mathcal{U}$ the problem
\begin{equation}\label{e2.7}
\left\{
\begin{array}{lll}
-\rho\Delta u=u-\left(u^{+}\right)^3\,\, &\text{in}\,\,  B(1+v),\\
u=0 &\text{on}\,\,  \partial B(1+v)
\end{array}
\right.\nonumber
\end{equation}
has a unique sign-changing solution $u=u_{\rho, v}\in C_{G}^{2,\alpha}(B(1+v))$ in a bounded neighborhood of $u_{\rho}$. Moreover the dependence of $u$ on the function $v$ is $C^1$ and $u_{\rho,0}=u_{\rho}$.}
\\

Define $F:\left(\rho_{0},\bar{\lambda}_{2}^{-1}\right)\times\mathcal{U}\rightarrow C_{G,m}^{1,\alpha}\left(\mathbb{S}^{N-1}\right)$ by
\begin{equation}\label{e2.8}
F(\rho,v)=\partial_{\nu}u|_{\partial B(1+v)}-\frac{1}{\text{Vol}(\partial B(1+v))}\int_{\partial B(1+v)}\partial_{\nu}u \,\text{d}S,\nonumber
\end{equation}
where $\mathcal{U}$ and $u=u_{\rho,v}$ are as given by Proposition 3.3. Notice that $F(\rho,v)=0$ if and only if $\partial_{\nu}u$ is constant on the boundary $\partial B(1+v)$ and $F\in C^1\left(\left(\rho_{0},\bar{\lambda}_{2}^{-1}\right)\times \mathcal{U}, C_{G,m}^{1,\alpha}\left(\mathbb{S}^{N-1}\right)\right)$. Obviously, $F(\rho,0)=0$ for all $\rho\in \left(\rho_{0},\bar{\lambda}_{2}^{-1}\right)$. In fact, we are aimed to find a branch of nontrivial solutions $(\rho,v)$ to the equation $F(\rho,v)=0$ bifurcating from some interval. Thus here we state an existence result as in Proposition 3.4 in \cite{Ruiz}.
\\ \\
\textbf{Lemma 3.4.} \emph{Assume that $\rho\in \left(\rho_{0},\bar{\lambda}_{2}^{-1}\right)$,
then for given $v\in  C^{2,\alpha}_G \left(\mathbb{S}^{N-1}\right)$, there exists a unique solution $\psi_{\rho,v}\in C_G^{2,\alpha}\left(B\right)$ of the problem}
\begin{equation}\label{e2.9}
\left\{
\begin{array}{ll}
-\rho\Delta\psi-\psi+3\left(u_{\rho}^{+}\right)^{2}\psi=0\,\, &\text{in}\,\, B,\\
\psi=v &\text{on}\,\, \partial B.
\end{array}
\right.
\end{equation}
\emph{In addition, if $v\in C^{2,\alpha}_G \left(\mathbb{S}^{N-1}\right)$, then $\psi_{\rho,v}\in E$ and $\int_{\partial B}\frac{\partial\psi_{\rho,v}}{\partial \nu} \,\text{d}S=0$.}
\\

Thus, for $\rho\in \left(\rho_{0},\bar{\lambda}_{2}^{-1}\right)$, we can define the linear continuous operator $H_{\rho}:C^{2,\alpha}_{G}\left(\mathbb{S}^{N-1}\right)\rightarrow C^{1,\alpha}_{G}\left(\mathbb{S}^{N-1}\right)$ by
\begin{equation}\label{e2.10}
H_{\rho}(v)=\partial_{\nu}\left(\psi_{\rho,v}\right)+(N-1)v,
\end{equation}
where $\psi_{\rho,v}$ is given by Lemma 3.4. In addition, from the Proposition 4.3 in \cite{Ruiz}, it's known that for any $\rho\in \left(\rho_{0},\bar{\lambda}_{2}^{-1}\right)$ and $v\in C_{G,m}^{2,\alpha}\left(\mathbb{S}^{N-1}\right)$, we have that
\begin{equation}\label{e2.11}
D_{v}(F)(\rho,0)=-\partial_{r}u(1)H_{\rho}(v),
\end{equation}
where $u=u_{\rho}$. Thus, a bifurcation of the branch $(\rho,0)$ of solutions to the equation $F(\rho,v)=0$ might appear only at some point $\left(\rho^*,0\right)$ such that $H_{\rho^*}$ becomes degenerate. To attain this aim, we first introduce a general function by
\begin{equation}
Q_{\rho}\left(\phi\right)=\int_{B}\left(\rho\vert\nabla \phi\vert^2-\phi^2+3\left(u_{\rho}^{+}\right)^2\phi^2\right)\,\text{d}x+(N-1)\rho\int_{\partial B}\phi^2\,\text{d}S.\nonumber
\end{equation}

By a perturbation argument as in \cite[Lemma 5.1 and Proposition 5.2]{Ruiz}, we have that
\\ \\
\textbf{Proposition 3.5.} \emph{For $\rho\in \left(\rho_{0},\bar{\lambda}_{2}^{-1}\right)$, there holds that}

(1) \emph{$\sigma_{1}(H_{\rho})=\min\left\{ \frac{1}{\rho}Q_{\rho}(\psi):\phi\in E, \int_{\partial B}\phi^2\,\text{d}S=1 \right\}$. Moreover, the infimum is attained};

(2) \emph{There exists some $\phi\in E$ and $\varepsilon>0$ such that $Q_{\rho}(\phi)<0$ for $\rho\in\left(\rho_{0},\rho_{0}+\varepsilon\right)$};

(3) \emph{There exists $\varepsilon>0$ such that for $\rho\in\left(\bar{\lambda}_{2}^{-1}-\varepsilon,\bar{\lambda}_{2}^{-1}\right)$, $Q_{\rho}|_{E}$ is positive definite}.
\\

Based on the Proposition 3.5, we have that $\sigma_{1}\left(H_{\rho}\right)$ can be negative for $\rho$ near $\rho_0$, but is positive for $\rho$ near $\bar{\lambda}_{2}^{-1}$. This provides a possibility for finding the zero point of $\sigma_{1}\left(H_{\rho}\right)$ with respect to $\rho$.

\section{Eigenvalues of the linearized operator}

\quad\,In this section, we will refine the eigenvalues $\sigma(\rho)$ of operator $H_{\rho}$ defined in (\ref{e2.10}) by using the spherical harmonics.
\\ \\
\textbf{Proposition 4.1.} \emph{The operator ${H}_{\rho}:C_{G,m}^{2,\alpha}\left(\mathbb{S}^{N-1}\right)\longrightarrow C_{G,m}^{1,\alpha}\left(\mathbb{S}^{N-1}\right)$ is a self-adjoint, first order elliptic operator.}
\\ \\
\textbf{Proof.} As defined in (\ref{e2.10}), the operator $H_{\rho}$ is the sum of the Dirichlet-to-Neumann operator for $-\rho\Delta-1+3\left(u_{\rho}^{+}\right)^{2}$ and
a constant times the identity, thus it is a first order elliptic operator.
On the other hand, let $\psi_1$ and $\psi_2$ be the solution of problem (\ref{e2.9})
with $v=v_1$ and $v=v_2$.
Then multiplying the equation of $\psi_1$ by $\psi_2$ and the equation of $\psi_2$ by $\psi_1$, integrating by parts, we obtain that
\begin{eqnarray}\label{e3.1}
\int_{\mathbb{S}^{N-1}}\left(\frac{\partial \psi_1}{\partial \nu}v_2-\frac{\partial \psi_2}{\partial \nu}v_1\right)\,\text{d}\theta=0.
\end{eqnarray}
Based on (\ref{e2.10}) and (\ref{e3.1}), we have that
\begin{eqnarray}
\int_{\mathbb{S}^{N-1}}\left(H_\rho\left(v_1\right)v_2-H_\rho\left(v_2\right)v_1\right)\,\text{d}\theta
=\int_{\mathbb{S}^{N-1}}\left(\frac{\partial \psi_1}{\partial \nu}v_2-\frac{\partial \psi_2}{\partial \nu}v_1\right)\,\text{d}\theta=0,\nonumber
\end{eqnarray}
which verifies the self-adjoint of $H_{\rho}$.\qed
\\

The necessary condition to bifurcate is that ${H}_{\rho}$ degenerates. In the following, we are devoted to finding the value of $\rho$ such that ${H}_{\rho}$ is degenerate. For any $v\in C_{G,m}^{2,\alpha}\left(\mathbb{S}^{N-1}\right)$, by virtue of the Fourier expansion with respect to spherical harmonics \cite[Theorem 3.2.11]{Groemer}, $v$ can be written as
\begin{equation}
v=\sum_{l=1}^\infty\sum_{j=1}^{m_l}a_{i_l,j}\zeta_{i_l,j}(\theta),\nonumber
\end{equation}
where $\zeta_{i_l,j}$ (normalized to $1$ in the $L^2$-norm) is an eigenfunction corresponding to $\mu_{i_l}$ with $\text{span}\left\{\zeta_{i_1,1}, \ldots,\zeta_{i_l,m_l}\right\}$ being the eigenspace. In the following, we will be concerned with the eigenvalue of ${H}_{\rho}$.
\\ \\
\textbf{Proposition 4.2.} \emph{For any $\rho>\rho_0$, ${H}_{\rho}$ possesses a strictly increasing sequence eigenvalues $\left\{\sigma_{i_l}(\rho)\right\}_{l\in \mathbb{N}}$ such that}
\begin{equation}
\sigma_{i_l}(\rho)=\sigma_{l}(H_{\rho}).\nonumber
\end{equation}
\emph{The eigenfunction corresponding to $\sigma_{i_l}(\lambda)$ is $\sum_{j=1}^{m_l}a_{i_l,j}\zeta_{i_l,j}(\theta)$ with $\sum_{j=1}^{m_l}a_{i_l,j}^2\neq0$.}
\\ \\
\textbf{Proof.} Let $\phi_0(r,\theta)=u_{\rho}'(r)v\left(\theta\right)$ with $r=\vert x\vert$,
it's easy to verify that
\begin{equation}
\Delta u_{\rho}'=\frac{N-1}{r^2}u_{\rho}'-\frac{u_{\rho}'-3\left(u_{\rho}^{+}\right)^{2}u_{\rho}'}{\rho}.\nonumber
\end{equation}
So we have that
\begin{eqnarray}
\Delta \phi_0&=&v\Delta u_{\rho}'+u_{\rho}'\Delta v\nonumber\\
&=&\left(\frac{N-1}{r^2}u_{\rho}'-\frac{u_{\rho}'-3\left(u_{\rho}^{+}\right)^{2}u_{\rho}'}{\rho}\right)v+
u_{\rho}'\Delta v\nonumber\\
&=&\sum_{l=1}^\infty\sum_{j=1}^{m_k}a_{i_l,j}\zeta_{i_l,j}(\theta)\left(\frac{N-1}{r^2}u_{\rho}'
-\frac{u_{\rho}'-3\left(u_{\rho}^{+}\right)^{2}u_{\rho}'}{\rho}\right)\nonumber\\
& &-\sum_{l=1}^\infty\sum_{j=1}^{m_l}a_{i_l,j}\mu_{i_l}\zeta_{i_l,j}(\theta)
u_{\rho}'\nonumber\\
&=&\sum_{l=1}^\infty\sum_{j=1}^{m_l}a_{i_l,j}u_{\rho}'\zeta_{i_l,j}(\theta)\left(\frac{N-1}{r^2}
-\frac{1-3\left(u_{\rho}^{+}\right)^{2}}{\rho}-\mu_{i_l}\right).\nonumber
\end{eqnarray}
It follows that
\begin{eqnarray}
-\rho \Delta \phi_0-\phi_0+3\left(u_{\rho}^{+}\right)^{2}\phi_0=\sum_{l=1}^\infty\sum_{j=1}^{m_l}a_{i_l,j}u_{\rho}'\zeta_{i_l,j}(\theta)\rho\left( -\frac{N-1}{r^2}+\mu_{i_1} \right).\nonumber
\end{eqnarray}

Let $\Psi$ be the unique solution of
\begin{eqnarray}
\left\{
\begin{array}{lll}
-\rho \Delta \Psi-\Psi+3\left(u_{\rho}\right)^{2}\Psi=\sum_{l=1}^\infty\sum_{j=1}^{m_l}a_{i_l,j}u_{\rho}'\zeta_{i_l,j}(\theta)\rho\left(\mu_{i_l}-\frac{N-1}{r^2}\right)\,\, &\text{in}\,\,  B,\\
\Psi=0  &\text{on}\,\, \partial B.
\end{array}
\right.\nonumber
\end{eqnarray}
Then we see that
\begin{equation}\label{phixt}
\Psi(r,\theta)=\phi_0(r,\theta)-\psi_v u_{\rho}'(1),
\end{equation}
where $\psi_v=\psi_{\rho,v}$ obtained in Lemma 3.4.
Moreover, we have that
\begin{equation}
\partial_r \Psi(1,\theta)=u_{\rho}''(1)v-u_{\rho}'(1) \partial_r\psi_v(1,\theta).\nonumber
\end{equation}
Note that $u_{\rho}''(1)+(N-1)u_{\rho}'(1)=0$. Thus we have that
\begin{equation}\label{e3.3}
\partial_r \Psi(1,\theta)=-(N-1)u_{\rho}'(1)v-u_{\rho}'(1) \partial_r\psi_v(1,\theta)=-u_{\rho}'(1)H_{\rho}(v).
\end{equation}
Let $V_l$ be the space spanned by the functions $\zeta_{i_1,j}(\theta)$, $\ldots$, $\zeta_{i_l,m_l}(\theta)$. We see that $H_{\rho}$ preserves $V_l$.
It follows that
\begin{equation}
{H}_{\rho} \sum_{j=1}^{m_l}a_{i_l,j}\zeta_{i_l,j}(\theta)=\sigma_{i_l}(\rho)\sum_{j=1}^{m_l}a_{i_l,j}\zeta_{i_l,j}(\theta),\nonumber
\end{equation}
where $\sigma_{i_l}(\rho)$ are the eigenvalues of $H_{\rho}$, and $\sum_{j=1}^{m_l}a_{i_l,j}\zeta_{i_l,j}(\theta)$ with $\sum_{j=1}^{m_l}a_{i_l,j}^2\neq0$ are the eigenfunctions associated to $\sigma_{i_l}(\rho)$.
So we have that
\begin{equation}
{H}_{\rho} \sum_{l=1}^n\sum_{j=1}^{m_l}a_{i_l,j}\zeta_{i_l,j}(\theta)=\sum_{l=1}^n\sigma_{i_l}(\rho)\sum_{j=1}^{m_l}a_{i_l,j}\zeta_{i_l,j}(\theta).\nonumber
\end{equation}
Since $H_{\rho}$ is bounded, we have that
\begin{equation}\label{ksaifenli1}
{H}_{\rho} v=\lim_{n\rightarrow+\infty}{H}_{\rho} \sum_{l=1}^{n}\sum_{j=1}^{m_l}a_{i_l,j}\zeta_{i_l,j}(\theta)=\sum_{l=1}^{+\infty}\sigma_{i_l}(\rho)\sum_{j=1}^{m_l}a_{i_l,j}\zeta_{i_l,j}(\theta).
\end{equation}
Hence the existence of $\sigma_{i_l}(\rho)$ has been proved. Next we study the properties of $\sigma_{i_l}(\rho)$.

With the help of (\ref{phixt}), we can write $\Psi(r,\theta)$ as
\begin{equation}
\Psi(r,\theta)=\sum_{l=1}^\infty\sum_{j=1}^{m_l}b_{i_l}(r)a_{i_l,j}\zeta_{i_l,j}(\theta),\nonumber
\end{equation}
where $b_{i_l}$ is the continuous solution on $(0,1]$ of
\begin{equation}
-\rho\left(\partial_r^2+\frac{N-1}{r}\partial_r\right)b+ \rho\mu_{i_l}b-\left(1-3\left(u_{\rho}\right)^{2}\right)b=
\rho u_{\rho}'(1)\left(\mu_{i_l}-\frac{N-1}{r^2}\right)\nonumber
\end{equation}
with $b_{i_l}(1)=0$.
In addition, from (\ref{e3.3}) we find that
\begin{equation}
-u_{\rho}'(1)\sigma_{i_l}(\rho)=\partial_r b_{i_l}(1)\nonumber
\end{equation}
for any $l\in \mathbb{N}$.

Combining (\ref{phixt}) with (\ref{ksaifenli1}), we have that
\begin{equation}
\psi_v=\sum_{l=1}^\infty\sum_{j=1}^{m_l}c_{i_l}(r)a_{i_l,j}\zeta_{i_l,j}(\theta),\nonumber
\end{equation}
where $c_{i_l}$ is the continuous solution on $(0,1]$ of
\begin{equation}\label{ckequation}
-\rho\left(\partial_r^2+\frac{N-1}{r}\partial_r\right)c+\rho\mu_{i_l}c-\left(1-3\left(u_{\rho}^+\right)^{2}\right)c=0
\end{equation}
with $c_{i_l}(1)=1$.
This implies that
\begin{equation}
\sigma_{i_l}(\rho)=c_{i_l}'(1)+N-1.\nonumber
\end{equation}

Now we mainly focus on the case of $l=1$ and let $\widetilde{c}(r,\theta)=c_{i_1}(r)\phi(\theta)$ with $\theta\in \mathbb{S}^{N-1}$, where
$\phi$ is a normalized eigenfunction corresponding to $\mu_{i_1}$.
So, using (\ref{ckequation}), we have that
\begin{equation}
\Delta \widetilde{c}=\phi(\theta)\Delta c_{i_1}+c_{i_1}\Delta \phi(\theta)=\frac{3\left(u_{\rho}^{+}\right)^2-1}{\rho}\widetilde{c}.\nonumber
\end{equation}
Multiplying the above equation by $\widetilde{c}$ and integrating by part, we get that
\begin{equation}
-\int_{\partial B}\widetilde{c} \nabla \widetilde{c}\cdot\nu\,\text{d}S+\int_{B} \left\vert\nabla \widetilde{c}\right\vert^2\,\text{d}x+\int_{ B}\frac{3\left(u_{\rho}^{+}\right)^2-1}{\rho}\widetilde{c}^2\,\text{d}x=0.\nonumber
\end{equation}
It follows that
\begin{equation}
c_{i_1}'(1)=\int_0^{1} r^{N-1}\left(c_{i_1}'^2+\frac{\left(3u_{\rho}^{+}\right)^2-1}{\rho}c_{i_1}^2\right)\,\text{d}r+\mu_{i_1}\int_0^{1} r^{N-3}c_{i_1}^2\,\text{d}r,\nonumber
\end{equation}
due to $c_{i_1}(1)=1$ and $c'_{i_1}(0)=0$.

For any $\psi\in E$ with $\int_{\partial B_1}\psi^2\,\text{d}S=1$,
there exist $\psi_l$ ($l\in \mathbb{N}\cup\{0\}$) and $a_{i_l,j}$ ($j\in\left\{1,\ldots,m_l\right\}$) with $\sum_{j=1}^{m_l}a_{i_l,j}^2=1$ for any $l$ such that
\begin{equation}
\psi(r,\theta)=\sum_{l=0}^{+\infty}\psi_l(r)\sum_{j=1}^{m_k}a_{i_l,j}\zeta_{i_l,j}(\theta).\nonumber
\end{equation}
When $l=0$, $\mu_{i_0}=0$ and $\zeta_{i_0,j}(\theta)$ is independent of $\theta$.
So, $\psi_{0}(1)=0$ which contradicts $\int_{\partial B_1}\psi^2\,dx=1$.
Then we have that
\begin{eqnarray}
\frac{Q_{\rho}}{\rho}\left(\psi\right)&=&\sum_{l=1}^{+\infty}\int_0^{1} r^{N-1}\left(\psi_{i_l}'^2-\frac{1-3(u_{\rho}^+)^{2}}{\rho}\psi_{i_l}^2\right)\,\text{d}r+\sum_{l=1}^{+\infty}\mu_{i_l}\int_0^{1} r^{N-3}\psi_{i_l}^2\,\text{d}r+N-1\nonumber\\
&\geq&\sum_{l=1}^{+\infty}\int_0^{1} r^{N-1}\left(\psi_{i_l}'^2-\frac{1-3\left(u_{\rho}^+\right)^{2}}{\rho}\psi_{i_l}^2\right)\,\text{d}r+\mu_{i_1}\sum_{l=1}^{+\infty}\int_0^{1} r^{N-3}\psi_{i_l}^2\,\text{d}r+N-1\nonumber\\
\nonumber\\
&\geq&\frac{1}{\rho}Q_{\rho}\left(\widetilde{\phi}\right),\nonumber
\end{eqnarray}
where
\begin{equation}
\widetilde{\phi}(r,\theta)=\sum_{l=1}^{+\infty}\psi_l(r)\sum_{j=1}^{m_1}a_{i_1,j}\zeta_{i_1,j}(\theta)\nonumber
\end{equation}
and $\int_{\partial B_1}\widetilde{\phi}^2\,dx=1$.
Therefore, the infimum of $\frac{Q_{\rho}(\phi)}{\rho}$ in $E$ with $\int_{\partial B}\phi^2\,dx=1$ is attained.
Define
\begin{equation}
\sigma_1\left(H_{\rho}\right):=\inf\left\{\frac{1}{\rho}Q_{\rho}(\phi):\phi\in E,\int_{\partial B}\phi^2\,\text{d}S=1\right\}.\nonumber
\end{equation}
We next to investigate the relations of $\sigma_1\left(H_{\rho}\right)$ and $\sigma_{i_1}(\rho)$.

The above argument implies that there exists $\phi(r,\theta)=\upsilon(r)\sum_{j=1}^{m_1}a_{i_1,j}\zeta_{i_1,j}(\theta)$ with $\sum_{j=1}^{m_1}a_{i_1,j}^2=1$ such that
$\phi\in E$, $\int_{\partial B}\phi^2\,\text{d}S=1$ and
\begin{equation}
\sigma_1\left(H_{\rho}\right)=\frac{1}{\rho}Q_{\rho}(\phi).\nonumber
\end{equation}
It follows that
\begin{equation}
\sigma_1\left(H_{\rho}\right)=\int_0^{1} r^{N-1}\left(\upsilon'^2-\frac{1-3\left(u_{\rho}^{+}\right)^{2}}{\rho}\upsilon^2\right)\,\text{d}r+\mu_{i_1}\int_0^{1} r^{N-3}\upsilon^2\,\text{d}r+N-1,\nonumber
\end{equation}
which is the functional of
\begin{equation}\label{weaksolution}
\left\{
\begin{array}{lll}
-\rho\left(\partial_r^2+\frac{N-1}{r}\partial_r\right)\upsilon+\frac{\rho\mu_{i_1}}{r^2(r)}\upsilon+\left(3\left(u_\rho^{+}\right)^2-1\right)\upsilon=0,\\
\upsilon(1)=1.
\end{array}
\right.
\end{equation}
So $\upsilon$ is a weak solution of (\ref{weaksolution}).
By Schauder elliptic estimates, $\upsilon$ is also the classical solution of (\ref{weaksolution}).
By Lemma 3.4 we deduce that $\upsilon(r)\equiv c_{i_1}(r)$.
So we get that
\begin{equation}
\sigma_1\left(H_{\rho}\right)=c_{i_1}'(1)+N-1=\sigma_{i_1}(\rho).\nonumber
\end{equation}
Thus the eigenspace corresponding to $\sigma_{i_1}(\rho)$ is just $V_1$. We next study the high eigenvalues.

For $\psi\in E$, we call $\psi\in E_k$ if there exist $\psi_{i_k}$ and $a_{i_k,j}$ such that
\begin{equation}
\psi(r,\theta)=\psi_{i_k}(r)\sum_{j=1}^{m_k}a_{i_k,j}\zeta_{i_k,j}(\theta).\nonumber
\end{equation}
Set
\begin{equation}
E_{k-1}^c:=E\setminus \cup_{i=1}^{k-1}E_{i}.\nonumber
\end{equation}
We will explore the infimum of ${Q}_\rho(\phi)/\rho$ in $E_{k-1}^c$ with $\int_{\partial B}\phi^2\,\text{d}S=1$.
For a function $\phi$ on $[1,+\infty)$ we define
\begin{equation}
Q_{\rho,k}(\phi)=\int_0^{1} r^{N-1}\left(\phi'^2-\frac{1-3\left(u_{\rho}^{+}\right)^{2}}{\rho}\phi^2\right)\,\text{d}r+\mu_{i_k}\int_0^{1} r^{N-3}\phi^2\,\text{d}r+N-1.\nonumber
\end{equation}
Since $\mu_{i_k}$ is increasing with respect to $k$, $Q_{\rho,k}(\phi)$ is increasing with respect to $k$.
For any $\psi\in E_k^c$ with $\int_{\partial B}\psi^2\,\text{d}S=1$,
there exist $\psi_l$ ($l\geq k$) and $a_{i_l,j}$ ($j\in\left\{1,\ldots,m_l\right\}$) with $\sum_{j=1}^{m_l}a_{i_l,j}^2=1$ for any $l\geq k$ such that
\begin{equation}
\psi(r,\theta)=\sum_{l=k}^{+\infty}\psi_l(r)\sum_{j=1}^{m_l}a_{i_l,j}\zeta_{i_l,j}(\theta).\nonumber
\end{equation}
Then we have that
\begin{eqnarray}
\frac{1}{\rho}{Q}_\rho\left(\psi\right)&=&\sum_{l=k}^{+\infty}\int_0^1 r^{N-1}\left(\psi_l'^2-\frac{1-3\left(u_{\rho}^{+}\right)^{2}}{\rho}\psi_l^2\right)\,\text{d}r+\sum_{l=k}^{+\infty}\mu_{i_l}\int_0^1 r^{N-3}\psi_{i_l}^2\,\text{d}r\nonumber\\
& &+(N-1)\nonumber\\
&\geq&\sum_{l=k}^{+\infty}\int_0^1 r^{N-1}\left(\psi_l'^2-\frac{1-3\left(u_{\rho}^{+}\right)^{2}}{\rho}\psi_l^2\right)\,\text{d}r+\mu_{i_k}\sum_{i=l}^{+\infty}\int_0^1 r^{N-3}\psi_l^2\,\text{d}r\nonumber\\
& &+(N-1)
\nonumber\\
&\geq&\frac{1}{\rho}{Q}_\rho\left(\widetilde{\phi}\right),\nonumber
\end{eqnarray}
where
\begin{equation}
\widetilde{\phi}(r,\theta)=\sum_{l=k}^{+\infty}\psi_l(r)\sum_{j=1}^{m_k}a_{i_k,j}\zeta_{i_k,j}(\theta)\in E_k\nonumber
\end{equation}
and $\int_{\partial B}\widetilde{\phi}^2\,\text{d}S=1$.
Therefore, the infimum is attained in $E_k$.
Hence we can set
\begin{equation}
\sigma_k\left(H_\rho\right):=\inf\left\{\frac{1}{\rho}{Q}_\rho(\phi):\phi\in E_{k-1}^c,\int_{\partial B}\phi^2\,\text{d}S=1\right\}.\nonumber
\end{equation}

There exists $\phi(r,\theta)=\upsilon(r)\sum_{j=1}^{m_k}a_{i_k,j}\zeta_{i_k,j}(\theta)$ with $\sum_{j=1}^{m_k}a_{i_k,j}^2=1$ such that
$\phi\in E_{k-1}^c$, $\int_{\partial B}\phi^2\,\text{d}S=1$ and
\begin{equation}
\sigma_k\left(H_\rho\right)=\frac{1}{\rho}{Q}_\rho(\phi).\nonumber
\end{equation}
It follows that
\begin{equation}
\sigma_k\left(H_\rho\right)=\int_0^1 r^{N-1}\left(\upsilon'^2-\frac{1-3\left(u_{\rho}^{+}\right)^{2}}{\rho}\upsilon^2\right)\,\text{d}r+\mu_{i_k}\int_0^1 r^{N-3}\upsilon^2\,\text{d}r+(N-1),\nonumber
\end{equation}
which is the functional of the following problem
\begin{equation}
\left\{
\begin{array}{lll}
-\rho\left(\partial_r^2+\frac{N-1}{r}\partial_r\right)\upsilon+\frac{\rho\mu_{i_k}}{r^2(r)}\upsilon+\left(3\left(u_\rho^{+}\right)^2-1\right)\upsilon=0,\\
\upsilon(1)=1.
\end{array}
\right.\nonumber
\end{equation}
By Lemma 3.4 we deduce that $\upsilon(r)\equiv c_{i_k}(r)$.
So we get that
\begin{equation}
\sigma_k\left(H_\lambda\right)=c_{i_k}'(1)+N-1=\sigma_{i_k}(\rho).\nonumber
\end{equation}
Since $E\varsupsetneq E_1^c\varsupsetneq \cdots\varsupsetneq E_k^c\varsupsetneq \cdots$ and $\sigma_k\left(H_\rho\right)$
is  attained in $E_k$, we have that
\begin{equation}
\sigma_1\left(H_\rho\right)<\sigma_2\left(H_\rho\right)<\cdots<\sigma_k\left(H_\rho\right)<\cdots.\nonumber
\end{equation}
It follows that
\begin{equation}
\sigma_{i_1}\left(\rho\right)<\sigma_{i_2}\left(\rho\right)<\cdots<\sigma_{i_k}\left(\rho\right)<\cdots,\nonumber
\end{equation}
which is the desired conclusion.
\qed
\\

We write $\sigma(\lambda)$ as $\sigma_{i_1}(\lambda)$ for simplicity.
Proposition 3.5 implies that $\sigma\left(\bar{\lambda}_2^{-1}\right)>0$ and $\sigma\left(\rho_0\right)<0$.
The definitions of ${Q}_\rho$ imply that $\sigma(\rho)$ is continuous.
So, $\sigma(\rho)$ has at least one zero.
We use $\rho_1$ to denote the biggest zero of $\sigma(\rho)$ such that $\sigma$ is negative in a small neighborhood on the left.
Clearly, one see that $\rho_1>\rho_0$. Taking $\rho_2>\rho_1$ such that $\sigma(\rho)>0$ for any $\rho>\rho_2$.
In view of Proposition 4.2, we have that
$\sigma_{i_k}\left(\rho\right)>0$ for any $k>1$ and $\rho>\rho_1$.

\section{Proof of Theorems 1.3}

\quad\,
We now show the desired conclusions of Theorem 1.3 by applying Theorem 1.2 or Corollary 2.1.
\\ \\
\textbf{Proof of Theorem 1.3.} Let $X=C_{G,m}^{2,\alpha}\left(\mathbb{S}^{N-1}\right)$ and $Y=C_{G,m}^{1,\alpha}\left(\mathbb{S}^{N-1}\right)$.
Then, clearly, the embedding of $X\hookrightarrow Y$ is compact.
Let $\mathcal{O}=\left(\rho_{0},\bar{\lambda}_{2}^{-1}\right)\times\mathcal{U}$.
We have known that $F:\mathcal{O}\rightarrow Y$ is $C^1$ with $F(\rho,0)=0$ for any $\rho\in \left(\rho_{0},\bar{\lambda}_{2}^{-1}\right)$.

We claim that $D_vF\left(\rho_1,0\right)$ is a Fredholm operator with index $0$. From \cite[Lemma 6.2]{Ruiz} we know that the kernel space of the
linearized operator $D_vF\left(\rho_1,0\right)$ is $V_1$.
So the kernel space of $D_vF\left(\rho_1,0\right)$ has odd multiplicity.
For any $v$ belonging to ${X}$ with $v>-1$, by virtue of the Fourier expansion with respect to spherical harmonics \cite[Theorem 3.2.11]{Groemer}, $v$ can be written as
\begin{equation}
v=\sum_{l=1}^\infty\sum_{j=1}^{m_l}a_{i_l,j}\zeta_{i_l,j}(\theta),\nonumber
\end{equation}
where $\zeta_{i_l,j}$ is an eigenfunction corresponding to $\mu_{i_l}$ with $\text{span}\left\{\zeta_{i_1,1}, \ldots,\zeta_{i_l,m_l}\right\}$ being the eigenspace.
By (\ref{ksaifenli1}), we know that
\begin{equation}
D_vF\left(\rho,0\right)v=\sum_{l=1}^{+\infty}\sigma_{i_l}(\rho)\sum_{j=1}^{m_l}a_{i_l,j}\zeta_{i_l,j}(\theta).\nonumber
\end{equation}
Since $\sigma\left(\rho_1\right)=0$ and $\sigma_{i_l}\left(\rho_1\right)>0$ for any $l>1$, we have
\begin{equation}
D_vF\left(\rho_1,0\right)v=\sum_{l=2}^{+\infty}\sigma_{i_l}(\rho)\sum_{j=1}^{m_l}a_{i_l,j}\zeta_{i_l,j}(\theta).\nonumber
\end{equation}
Consequently, the image of $D_vF\left(\rho_1,0\right)$ is the closure of
\begin{equation}
\bigoplus_{i\geq2}V_i.\nonumber
\end{equation}
Therefore, we obtain that
\begin{equation}
\text{dim} \text{Ker}\left(D_vF\left(\rho_1,0\right)\right)=\text{codim} \text{Im}\left(D_vF\left(\rho_1,0\right)\right)=m_1.\nonumber
\end{equation}
So, to verify this claim, it is enough to show that $D_vF\left(\rho_1,0\right)$ has a closed range.

For $\psi_{\rho,v}$ is given by Lemma 3.4, by \cite[Theorem 8.16]{Gilbarg}, we have that $\sup_{B}\left\vert \psi_{\rho,v}\right\vert\leq \sup_{\partial B}\vert v\vert$.
Furthermore, applying \cite[Theorem 6.6]{Gilbarg}, we have that
\begin{equation}
\Vert \psi_{\rho,v}\Vert_{C_G^{2,\alpha}\left(B\right)}\leq  M \Vert v\Vert_{C_{G,m}^{2,\alpha}\left(\mathbb{S}^{N-1}\right)}
\nonumber
\end{equation}
for some $M>0$.
It follows that $\left\Vert {H}_\rho(v)\right\Vert_{C_{G,m}^{1,\alpha}\left(\mathbb{S}^{N-1}\right)}\leq M_1 \Vert v\Vert_{C_{G,m}^{2,\alpha}\left(\mathbb{S}^{N-1}\right)}$ for some $M_1>0$.
Thus, the operator
\begin{equation}
{H}_\rho:C_{G,m}^{2,\alpha}\left(\mathbb{S}^{N-1}\right)\longrightarrow C_{G,m}^{1,\alpha}\left(\mathbb{S}^{N-1}\right)\nonumber
\end{equation}
is bounded.
Assume that $v_n\in X$ and $y_n\in Y$ with $D_vF\left(\rho_1,0\right)v_n=y_n\rightarrow y_0$ in $Y$.
We sue $P$ to denote the projection of $X$ into $V_1$.
Then we have that $D_vF\left(\rho_1,0\right)(I-P)v_n=y_n$.
It follows that
$(I-P)v_n=\left(D_vF\left(\rho_1,0\right)\right)^{-1}y_n$.
By Banach inverse operator theorem, $\left(D_vF\left(\rho_1,0\right)\right)^{-1}:R\left(\left(D_vF\left(\rho_1,0\right)\right)^{-1}\right)\rightarrow (I-P)X$ is bounded.
Thus, $(I-P)v_n$ is convergent. Since $Pv_n\in V_1$ and is bounded, it is also convergent.
So we have that $D_vF\left(\rho_1,0\right)v_0=y_0$ for some $v_0\in X$. Hence, $D_vF\left(\rho_1,0\right)$ has a closed range.
Similarly, we also have that $D_vF\left(\rho,0\right)$ has a closed range for any $\rho\in\left(\rho_{0},\bar{\lambda}_{2}^{-1}\right)$.
Since $F$ is $C^1$, by \cite[Theorem 5.17]{Kato}, $D_vF\left(\rho,0\right)$ is also a Fredholm operator with index $0$.

More general, since $F$ is $C^1$ and $\mathcal{O}$ is bounded, $D_vF\left(\rho,v\right)$ is bounded for any $\left(\rho,v\right)\in\mathcal{O}$.
Then, reasoning as the above, for any $\left(\rho,v\right)\in\mathcal{O}$, we also can show that $D_vF\left(\rho,v\right)$ is also a Fredholm operator with index $0$.
By Proposition 2.1 $F$ is locally proper.
Since $\sigma(\rho)>0$ for any $\rho>\rho_2$, we see that the $0$-group index
\begin{equation}
\Sigma\left(\rho_2+\varepsilon\right)=1\nonumber
\end{equation}
for any $\varepsilon>0$ small enough.
From the argument of \cite[Theorem 2.1]{Ruiz} we know that the $0$-group index
\begin{equation}
\Sigma\left(\rho_1-\varepsilon\right)=(-1)^{m_1}=-1.\nonumber
\end{equation}

Applying Corollary 2.1, it can be obtained that
$\mathcal{S}$
possesses a component $\mathscr{C}_{1}$ emanating from $\left[\rho_1,\rho_2\right]\times \{0\}$, such that either $\mathscr{C}_{1}$ contains a point $\left(\mu,0\right)\in \mathcal{O}$ with $\mu\not\in\left[\rho_1,\rho_2\right]$ or $\mathscr{C}_{1}\cap \partial\mathcal{O}\neq\emptyset$, \noindent where $\mathcal{S}$ denotes the closure of the set of nontrivial solution pairs of $F(\rho,v)=0$ in $\mathcal{O}$.
 \qed


\bibliographystyle{amsplain}
\makeatletter
\def\@biblabel#1{#1.~}
\makeatother


\providecommand{\bysame}{\leavevmode\hbox to3em{\hrulefill}\thinspace}
\providecommand{\MR}{\relax\ifhmode\unskip\space\fi MR }
\providecommand{\MRhref}[2]{%
  \href{http://www.ams.org/mathscinet-getitem?mr=#1}{#2}
}
\providecommand{\href}[2]{#2}

\end{document}